\documentclass[a4paper,11pt]{article}
\usepackage{amssymb,latexsym}
\newtheorem{thm}{Theorem}[section]

\newtheorem{cor}[thm]{Corollary}
\newtheorem{prop}[thm]{Proposition}
\newtheorem{conj}[thm]{Conjecture}
\newtheorem{rem}[thm]{Remark}

\input epsf

\title{Twisting the Stern sequence}
\author{Roland Bacher}

\begin{document}
\maketitle




\begin{abstract} We describe a few features of the Stern sequence and 
of a closely related sequence obtained by adding a sign-twist in 
the recursive definition
of the Stern sequence.\footnote{Keywords: 
Stern sequence, automatic sequence, regular sequence. 
Math. class: 11B85}
\end{abstract}

\section{Main results}

In this paper, we identify a complex 
sequence $(a_n)_{n\in \mathbb N}$ with the corresponding function
$a:\mathbb N\longmapsto \mathbb C$. We write thus always $a(n)$ instead of
$a_n$.

The Stern sequence or Stern-Brocot sequence with first terms given by 
$$0,1,1,2,1,3,2,3,1,4,3,5,2,5,3,4,1,5,4,7,3,8,5,7,2,7,5,8,3,\dots$$
(cf. sequence A2487 of \cite{OEIS})
is the integral sequence $s:\mathbb N\longrightarrow \mathbb N$
recursively defined by $s(0)=0,s(1)=1$ and
$s(2n)=s(n),s(2n+1)=s(n)+s(n+1)$ for $n\geq 1$.
It is closely related to the Farey tree and induces a one-to-one map 
$n\longmapsto s(n)/s(n+1)$ between $\mathbb N$ and
non-negative rational numbers, cf. \cite{CW} or
Chapter 16 of \cite{AZ}.
It is also an example of a $2-$regular sequence, see Chapter 16 of 
\cite{AS}. 

The following result gives a different, perhaps not very well-known,
description of the Stern sequence.

\begin{prop}\label{propotherdescr}
$s(n)$ equals the number of distinct subsequences of the form 
$1,101,10101,\dots=\{1(01)^*\}$ in the binary expansion 
$\epsilon_l\dots\epsilon_1\epsilon_0$ of $n=\sum_{k=0}^l\epsilon_l2^k$
(where $\epsilon_0,\dots,\epsilon_l\in\{0,1\}$).
\end{prop}

Proposition \ref{propotherdescr} is in fact a particular case of 
Proposition \ref{propnoncommprod}, an easy result concerning
rational series in non-commuting variables.

{\bf Example} The binary expansion $1011$ of $11=2^3+2^1+2^0$ 
contains the following five subsequences (highlighted by bold letters)
$${\mathbf 1}\mathit{011},\mathit{10}{\mathbf 1}\mathit{1},\mathit{101}
{\mathbf 1},\mathbf{101}\mathit{1},
\mathbf{10}\mathit{1}{\mathbf 1}$$
of the form $1(01)^*$. We have 
$s(11)=s(5)+s(6)=s(2)+s(3)+s(3)=s(1)+2(s(1)+s(2))=5s(1)=
5$.

Proposition \ref{propotherdescr} allows to parametrize Stern
sequences by counting a subsequence of the form 
$1(01)^k$ with weight $w^k$, see Proposition \ref{propparamstern}
for formulae. For $n=11$ we get for instance $3+2w$.

In this paper we introduce a related sequence which will be 
called the {\it twisted Stern
sequence} since it is obtained by twisting the 
recursive definition of the Stern sequence
with a sign. More precisely, we define 
the twisted Stern sequence $t(0),t(1),\dots$ recursively by 
$t(0)=0,t(1)=1$ and
$t(2n)=-t(n),t(2n+1)=-t(n)-t(n+1)$ for $n\geq 1$. It starts as
$$0,1,-1,0,1,1,0,-1,-1,-2,-1,-1,0,1,1,2,1,3,2,3,1,2,1,1,0,-1,\dots\ .$$
An inspection of these first 
few terms shows already some striking similarities between the Stern 
sequence and its twisted relative. 
The aim of this paper is to describe a few properties of the Stern sequence 
and its twist. 

The following result (the identity for $s(n)$ in
assertion (i) is probably well-known to the experts)
is an illustration of the similarities between these two sequences:

\begin{thm} \label{thmstid} (i) We have 
\begin{eqnarray*}
s(2^e+n)&=&s(2^e-n)+s(n)\\
t(2^e+n)&=&(-1)^e\left(s(2^e-n)-s(n)\right)
\end{eqnarray*}
for all $e\geq 0$ and for all $n$ such that $0\leq n\leq 2^e$.

(ii) We have
$$t(3\cdot 2^e+n)=t(6\cdot 2^e-n)=(-1)^e s(n)$$ 
for all $e\geq 0$ and for all $n$ such that 
$0\leq n\leq 2^{e+1}$.
\end{thm}

The failure for $n>2^{e+1}$ of the formula
$$t(3\cdot 2^e+n)=(-1)^e s(n)$$
given by assertion (ii) can perhaps be
mended by the following conjectural identity based on
experimental observations.

\begin{conj}\label{conjgen}
There exists an integral sequence $u(0),u(1),u(2),\dots$ such that
we have
$$\sum_{n=0}^\infty t(3\cdot 2^e+n)z^n=(-1)^e
\left(\sum_{n=0}^\infty u(n)z^{n\cdot 2^e}\right)
\left(\sum_{m=0}^\infty s(m) z^m\right)$$
for all $e\in \mathbb N$.
\end{conj}

If the conjecture holds, the ordinary generating function of the
sequence $u(0),u(1),\dots$ is given by
$$\sum_{n=0}^\infty u(n)z^n=\frac{\sum_{n=0}^\infty
t(3+n)z^n}
{\sum_{n=0}^\infty s(n)z^n}$$
and it starts as 
$$1-2z^2-2z^5+4z^6+2z^7-6z^8+4z^9+2z^{10}-6z^{11}+8z^{12}+\dots$$

The first equality in assertion (ii) of Theorem \ref{thmstid}
shows that the finite
sequences $(-1)^et(3\cdot 2^e),(-1)^et(3\cdot 2^e+1,\dots,(-1)^e
t(6\cdot 2^e)$ of length 
$3^e+1$ are palindromic sequences of natural integers.
The first few such sequences are
$$\begin{array}{cccccccccccccccccccccccccccccccccccccccc}
0&&&&&&&&{\mathbf 1}&&&&&&&&{\mathbf 1}&&&&&&&&0\\
0&&&&1&&&&{\mathbf 1}&&&&2&&&&{\mathbf 1}&&&&1&&&&0\\
0&&1&&1&&2&&{\mathbf 1}&&3&&2&&3&&{\mathbf 1}&&2&&1&&1&&0\\
0&1&1&2&1&3&2&3&{\mathbf 1}&4&3&5&2&5&3&4&{\mathbf 1}&
3&2&3&1&2&1&1&0\end{array}
$$
with boldfaced $1$'s at one third and two thirds 
highlighting the underlying partial self-similarity structure.
The subsequence lying between the two boldfaced $1'$s 
appears also at the beginning of \cite{L}.  
One notices that all sequences start and end with zero and that all 
existing central elements are equal to $2$.

The polynomials defined by these palindromic sequences
are described by the following result:

\begin{thm} \label{thmtpoly} The polynomials 
$$\psi_e=(-1)^e\sum_{n=0}^{3\cdot 2^e}t(3\cdot 2^e+n)z^n$$
have the factorisations
\begin{eqnarray}
\psi_e&=&z(1+z^{2^e})(1+z+z^2)^e\prod_{n=0}^{e-2}
(1-z^{2^n}+z^{2^{n+1}})^{e-1-n}\\
&=&z(1+z^{2^e})\prod_{n=0}^{e-1}(1+z^{2^n}+z^{2^{n+1}})\ .\label{thmsecfact}
\end{eqnarray}
\end{thm}

\begin{rem}\label{remevalsert}
Theorem \ref{thmtpoly} implies the identity
$$\sum_{n=0}^\infty t(n)z^n=z-z^2+\sum_{e=0}^\infty (-1)^ez^{3\cdot 2^e+1}
(1+z^{2^e})\prod_{n=0}^{e-1}(1+z^{2^n}+z^{2^{n+1}})\ .$$
\end{rem}

Assertion (ii) of Theorem \ref{thmstid} yields 
$\lim_{e\rightarrow\infty} \psi_e=\sum_{n=0}^\infty s(n)t^n$. 
The factorisation (\ref{thmsecfact}) of Theorem \ref{thmtpoly}
gives a new proof of the following result due to Carlitz (see \cite{C}):

\begin{cor} \label{corCarl}
We have $\sum_{n=0}^\infty s(n)z^n=z\prod_{n=0}^\infty (1+z^{2^n}+
z^{2^{n+1}})$.
\end{cor}

(A direct proof of Corollary \ref{corCarl} is straightforward: 
The series $U(z)$ defined by the right-hand-side
starts as $z+\dots=s(0)+s(1)z+\dots$ and its even, respectively odd,
subseries are given by $U(z^2)$, respectively $\left(\frac{1}{z}+z\right)
U(z^2)$. Its coefficients satisfy thus the same recursion relations
as the elements of the Stern sequence.)

The Carlitz factorisation of Corollary \ref{corCarl} 
implies that $\sum_{n=0}^\infty s(n)z^n$ has no non-zero roots
in the open unit disc. This is not true for the ordinary generating series
$\sum_{n=0}^\infty t(n)z^n$ of the twisted Stern sequence which has
(infinitely?) many non-zero roots in the open unit disc.

Given a natural integer $k\geq 2$ and a natural integer $i$, 
we consider the endomorphism $\rho(i)$ of the vector-space (or module) 
of formal power series defined by 
$$\rho(i)\left(\sum_{n=0}^\infty a(n)z^n\right)=
\sum_{n=0}^\infty a(i+nk)z^n\ .$$
The $k-$kernel of a formal power series $A$ is the smallest vector space
(or module when working over a ring) $\mathcal V$ 
containing $A$ such that $\rho(0)\mathcal V,\dots,\rho(k-1)\mathcal V
\subset\mathcal V$.
A formal power series is $k-$regular if its $k-$kernel is finitely 
generated. Easy examples of $k-$regular series
are polynomials and ordinary generating series of periodic 
sequences. $k-$regular power series form a vector space (or module)
which is preserved by many natural operations such as
derivation, product, Hadamard product, shuffle product, ..., see 
\cite{AS} for details. The set of $k-$regular sequences with 
coefficients contained in a 
finite set (eg. in a finite field) coincides with the set of so-called
$k-$automatic sequences, see Theorem 16.1.5 of \cite{AS}.

A sequence $a(0),a(1),\dots $ is called $k-$regular if its
ordinary generating series $\sum_{n=0}^\infty a(n)z^n$ is a $k-$regular 
formal power series.

The following result is also a consequence of the Carlitz factorisation:

\begin{thm} \label{thmlogder}
The logarithmic derivative 
$$H(z)=\frac{d}{dz}\mathop{log}\left(\sum_{n=0}^\infty s(n+1)z^n\right)=
\frac{\sum_{n=1}^\infty ns(n+1)z^{n-1}}{\sum_{n=0}^\infty s(n+1)z^n}$$
of $\sum_{n=0}^\infty 
s(n+1)z^n$ is $2-$regular. More precisely, $H(z)$ is defined by 
the functional equation
$$H(z)=\frac{1+2z}{1+z+z^2}+2zH(z^2)\ .$$
\end{thm}

Coefficients of $H(z)$ appear as sequence A163659 in \cite{OEIS}.

$2-$regularity of the logarithmic derivation $H(z)$ is a special case
of the following result, perhaps already known to Sch\"utzenberger:

\begin{thm} \label{thmaffpres} Given $d$ $k-$regular series $A_1(z),\dots,A_d(z)$ over
some commutative ring $R$, $d$ linear forms 
$L_1(x_1,\dots,x_d),\dots,L_d(x_1,\dots,x_d)$ in $d$ unknowns with
coefficients in $R[z]$ and $d$ constants $\alpha_1,\dots,\alpha_d\in R$
such that $\alpha_i=A_i(0)+L_i(\alpha_1,\dots,\alpha_d)\pmod z$ 
for $i=1,\dots,d$, the system of equations
$$\begin{array}{l}
U_1(z)=A_1(z)+L_1(U_1(z^k),\dots,U_d(z^k)),\\
\qquad \vdots\\
U_d(z)=A_n(z)+L_n(U_1(z^k),\dots,U_d(z^k))\end{array}$$
determines a unique set of $d$ $k-$regular sequences
$U_1(z),\dots,U_d(z)$ with constant coefficients $\alpha_i=U_i(0)$
for $i=1,\dots,d$. 
\end{thm}

\begin{rem} We have $\alpha_i=A_i(0)$ if the linear form $L_i$
has all its coefficients in $zR[z]$.
\end{rem}

The series
$$A(z)=\prod_{n=0}^\infty\frac{1}{1-z^{2^n}}$$
satisfying $A(z)=\left(\sum_{n=0}^\infty z^n\right)A(z^2)$
is not $2-$regular (see Remark \ref{remnot2reg} below). This shows
that Theorem \ref{thmaffpres} can not be extended to
equations with linear forms having 
$k-$regular series as coefficients.

Theorem \ref{thmlogder} can be generalised as follows:

\begin{thm} \label{thmPreg} 
Let $P(z)$ be a polynomial with constant coefficient 
$1$. Then the series
$$A=\prod_{n=0}^\infty P(z^{k^n})$$
is $k-$regular.

Moreover, if all roots of $P(z)$ are complex roots of $1$ having finite
order, then the logarithmic derivative $B=A'/A$ of $A$ is also $k-$regular.
\end{thm}

\begin{rem} \label{remnot2reg}
Given a $k-$regular series $A(z)$ with constant coefficient $1$, 
the product $\prod_{n=0}^\infty A(z^{k^n})$ is generally not $k-$regular. 
Indeed, starting with the $2-$regular series 
$A(z)=1+z+z^2+\dots=\frac{1}{1-z}$,
the coefficient $b(n)$ in the series 
$$B=\sum_{n=0}^\infty b(n)z^n=\prod_{n=0}^\infty \frac{1}{1-z^{2^n}}$$
counts the number of partitions of $n$ into powers of $2$, see 
sequence A123 in \cite{OEIS}, and 
$\mathop{log}(b(2n))$ is asymptotically equal to 
$\frac{1}{2\mathop{log}2}\left(\mathop{log}\frac{n}{\mathop{log}n}\right)^2$
(see equation 1.3 in \cite{dB}) which is 
incompatible with $k-$regularity of $B$ by Theorem 16.3.1
in \cite{AS}.
\end{rem}

\begin{rem} Another famous sequence illustrating 
Theorem \ref{thmPreg} is the sequence
$$\prod_{n=0}^\infty (1-z^{2^n})=\sum_{n=0}^\infty (-1)^{tm(n)}z^n$$
related to the Thue-Morse sequence $n\longmapsto tm(n)$ defined by 
digit-sums modulo $2$ for binary expansions of natural integers.
\end{rem}

Another link between the two sequences $s$ and $t$ 
is given by determinants of $2\times 2-$matrices. 
For $n\geq 1$ we consider the matrix
$$M(n)=\left(\begin{array}{cc}s(n)&s(n+1)\\t(n)&t(n+1)\end{array}\right)$$
with first row two consecutive terms of $s$ and and second row 
the two corresponding consecutive terms of $t$. 
The first matrices are
$2\times 2-$submatrices defined by two consecutive rows of 
$$\begin{array}{rrrrrrrrrrrrrrrrrrrrrrrrrrrrrrrrrrrr}
1&1&2&1&3&2&3&1&4&3&5&2&5&3&5&1&5&4&7\\
1&-1&0&1&1&0&-1&-1&-2&-1&-1&0&1&1&2&1&3&2&3
\end{array}$$

Injectivity of the map $n\longmapsto s(n+1)/s(n)$ shows that the matrices
$M(n)$ are all different. Their determinants 
are characterised by the following result.

\begin{thm} \label{thmmatr} 
We have $\vert \det(M(n))\vert =2$ for all $n\geq 1$.
More precisely,
$$\det(M(n))=-2(-1)^k\hbox{ if }2^k\leq n<2^{k+1}\ .$$
In particular, all matrices $M(n)$ are invertible for $n\geq 1$.
\end{thm}

Since the recursive definitions of the integral sequences $s$ and $t$ 
differ only by signs (and since they satisfy the same initial conditions)
they have the same reduction modulo $2$ characterised by the following 
(easy) result, already contained in \cite{S}:

\begin{prop} \label{propstmod2}
The integers $s(n)$ and $t(n)$ are even if and only if
$n$ is divisible by $3$.
\end{prop}

The determinants of the non-singular matrices $M(n)$ are thus in some 
sense as small as
possible: Indeed, since the sequences $s$ and $t$ coincide modulo $2$,
a matrix $M(n)$ involves either a column consisting of even integers
or all its four entries are odd integers and such matrices
have even determinants.

We end this paper with a last result, going back to Stern (see \cite{S})
for the first part of assertion (i):

\begin{thm}\label{thmdiv} (i) The integer $s(n)$ divides $s(n-1)+s(n+1)$ for every 
$n\geq 1$. More precisely, we have 
$$\frac{s(n-1)+s(n+1)}{s(n)}=1+2v_2(n)$$
where the $2-$valuation $v_2(n)$ is defined as the exponent of the highest 
power of $2$ dividing $n$.

The function 
$$C(z)=\sum_{n=1}^\infty \frac{s(n-1)+s(n+1)}{s(n)}z^n$$
satisfies $C(0)=0$ and 
$$C(z)=z\frac{1+2z}{1-z^2}+C(z^2)$$
and is $2-$regular.

\ \ (ii) The integers $t(n)$ and $t(n-1)+t(n+1)$ are both zero if 
$n=3\cdot 2^k$. They are both non-zero otherwise and $t(n)$ divides 
$t(n-1)+t(n+1)$. More precisely, we have
$$\frac{t(n-1)+t(n+1)}{t(n)}=1+2v_2(n)$$
if $n\not\in \{2^{\mathbb N},3\cdot 2^{\mathbb N}\}$,
$(t(0)+t(2))/t(1)=-1$ and
$$\frac{t(2^e-1)+t(2^e+1)}{t(2^e)}=1+2(e-2)$$
for all $e$ such that $e\geq 1$.
\end{thm}

The rest of this paper is organized as follows: The next section 
contains a proof of Proposition \ref{propotherdescr} and a few complements.

Section \ref{sectproofsmain} is devoted to
the (easy) proofs of Theorem \ref{thmstid}
and to a few more formulae and conjectures involving the 
Stern sequence and its twist. 

Section \ref{sectproofkreg} contains the proofs of Theorems
\ref{thmtpoly}, \ref{thmlogder}, \ref{thmaffpres}, \ref{thmPreg}
and Corollary \ref{corCarl}.

Section \ref{sectproofmatr} contains the easy proof of Theorem \ref{thmmatr}
and a few related results.

Section \ref{sectdiv} consists of  the short proof  
of Theorem \ref{thmdiv}.


\section{Proof of Proposition \ref{propotherdescr} and a few comments}

We give first a bijective proof of Proposition \ref{propotherdescr}.
We describe then briefly a weighted version of the Stern sequence
counting subsequences of the form $1(01)^*$ with weights encoding their
length. We give also a generalisation of Proposition \ref{propotherdescr}
therefore providing a (sketch of a) 
second proof for Proposition \ref{propotherdescr}.

{\bf Proof of Proposition \ref{propotherdescr}} 
We call a subsequence of a binary expansion $B(n)$ {\it admissible} if it 
is of the form $1(01)^*$. For example, 
$B(11)=1011=\epsilon_3\epsilon_2\epsilon_1
\epsilon_0$ has five admissible subsequences 
given by the set 
$$\{\epsilon_3,\epsilon_1,\epsilon_0,\epsilon_3\epsilon_2\epsilon_1,
\epsilon_3\epsilon_2\epsilon_0\}\ .$$
Since the number $b(n)$ of such subsequences satisfies
clearly $b(2n)=b(n)$, the equality $s(2n)=s(n)$ shows that 
we can restrict our attention to $n$ odd. We 
consider the two cases $n=4n+1$ and $n=4n-1$.

If $w$ is an admissible subsequence of $B(4n+1)$, then
the digit $\epsilon_1=0$ of the binary expansion of 
$4n+1$ is either contained in $w$ or not.
In the first case, admissibility of $w$ shows that $w$ 
contains also the last digit $\epsilon_0$ of $B(4n+1)$ and removal of 
$\epsilon_1\epsilon_0=01$ from $w$ yields a bijection between such 
admissible subsequences and 
admissible subsequences of $B(2n)$. In the second case
where $\epsilon_1$ is not involved in $w$,
we get a bijection between such  admissible subsequences and 
admissible subsequences
of the binary expansion $B(2n+1)=\dots\epsilon_3\epsilon_2\epsilon_0$ 
of the integer $2n+1$. This shows the identity $b(4n+1)=
b(2n)+b(2n+1)$.

We consider finally the case of an admissible subsequence $w$ of $B(4n-1)$.
If $\epsilon_0$ is not contained in $w$, then $w$ can be associated
with an admissible subsequence of $B(4n-2)$ or equivalenty of $B(2n-1)$.
Denoting by $l$ the least integer such that $B(4n-1)=\alpha 01^l$,
we consider now an admissible subsequence $w$ of $B(4n-1)$
which contains $\epsilon_0$. If the admissible subsequence
$w$ is not of the form $\beta \epsilon_l\epsilon_0$, we transform
it into the 
admissible subsequence $\beta\tilde{\epsilon}_l=\beta 1$ of $B(4n)=
\alpha\tilde{\epsilon}_l\tilde{\epsilon}_{l-1}\dots\tilde{\epsilon}_0=
\alpha 10^l$
or equivalently of $B(2n)=\alpha 10^{l-1}$ (obtained from $B(4n)$ 
by erasing the last 
digit $0$ never involved in an admissible subsequence). 
If $w=\beta \epsilon_l\epsilon_0$
we transform it into the admissible subsequence $\beta$ of $B(4n)$ or 
equivalently of $B(2n)$. 
This shows $b(4n-1)=b(2n-1)+b(2n)$ and ends 
the proof. \hfill$\Box$

\subsection{A weighted variation of the Stern sequence}

We denote by $S(n)\in\mathbb N[w]$ the weighted number of 
subsequences of the form 
$1(01)^*$ in the binary expansion $B(n)$ of $n$, giving the weigth
$w^k$ to a subsequence of the form $1(01)^k$. 
Similarly, we introduce 
$S_e(n)\in\mathbb N[w]$ as the weighted number of 
subsequences of the form 
$(10)^*$ in the binary expansion $B(n)$ of $n$, with weigth
$w^k$ for a subsequence of the form $(10)^k$.

\begin{prop} \label{propparamstern}
(i) Evaluating the polynomial $S(n)\in\mathbb N[w]$ at $w=1$ 
yields the Stern sequence.

\ \ (ii) The sequences $S(n)$ and $S_e(n)$ are uniquely determined by
the initial conditions $S(0)=0,S(1)=S_e(0)=S_e(1)=1$
and the recursive formulae
\begin{eqnarray*}
S(2n)&=&S(n)\\
S(2n+1)&=&S(n)+S_e(n)\\
S_e(2n)&=&wS(n)+S_e(n)\\
S_e(2n+1)&=&S_e(n)\ .
\end{eqnarray*}

\ \ (iii) The sequence $S(n)$ is also uniquely determined by the initial 
conditions 
$S(0)=0,S(1)=1$ and by the recursive formulae
\begin{eqnarray*}
S(2n)&=&S(n)\\
S(4n+1)&=&wS(2n)+S(2n+1)\\
S(4n-1)&=&S(2n-1)+S(2m+1)+(w-1)S(2m)
\end{eqnarray*}
where $n=2^a(2m+1)$.
\end{prop}

\begin{rem} Klavzar, Milutinovic and Petri have studied a different 
family of polynomials closely related to the Stern sequence by 
considering $B_0=0, B_1=1, B_{2n}=tB(n)$ and $B_{2n+1}=B_n+B_{n+1}$,
see \cite{KMP} for details.
\end{rem}

{\bf Proof of Proposition \ref{propparamstern}} Assertion (i) is obvious.

In the sequel, we use the notation introduced above during the proof of 
Proposition \ref{propotherdescr}.

The initial values for $S(n)$ and $S_e(n)$ in assertion (ii) are easy to 
check. The identity $S(2n)=S(n)$ is obvious since admissible 
subsequences of $B(2n)$ never involve the last digit $\epsilon_0=0$
in the binary expansion $B(2n)$ of $2n$.

Admissible subsequences of $B(2n+1)$ not containing the 
last digit $\epsilon_0$ of $B(2n+1)$ are in weight-preserving
bijection with admissible subsequences of $B(2n)$ or of $B(n)$. 
Removal of $\epsilon_0$ induces a weight-preserving bijection between
admissible subsequences of $B(2n+1)$ involving the last digit $
\epsilon_0$ of $B(2n+1)$ and
monomial contributions to $S_e(n)$. This proves $S(2n+1)=S(n)+S_e(n)$.

Monomial contributions to $S_e(2n)$ not involving the last digit
$\epsilon_0$ of $B(2n)$ are in (weight-preserving) bijection with 
monomial contributions to $S_e(n)$. Removing the last digit
of monomial contributions to $S_e(2n)$ involving the last digit
$\epsilon_0$ of $B(2n)$ yields admissible subsequences of $B(n)$
with weight reduced by $1$. This shows $S_e(2n)=S_e(n)+wS_e(n)$.

The identity $S_e(2n+1)=S_e(n)$ is due to the fact that monomial
contributions to $S_e(2n+1)$ never involve the last digit $\epsilon_0=1$
of $B(2n+1)$.

Assertion (iii) follows from the bijections used in the proof of 
Proposition \ref{propotherdescr}. We leave the details to the reader.
\hfill$\Box$

\subsection{Counting weighted subsequences and subfactors}

A famous result by Sch\"utzenberger implies essentially an
identification of $k-$regular sequences with the set of rational
formal power series in $k$ non-commuting variables. (One has to
be a little careful with leading zeros. A way of dealing with them
is to consider only formal power series involving no monomials
starting with the variable $x_0$ associated to the digit $0$.)

Proposition \ref{propotherdescr} is then a particular case
of the following well-known result which we give without proof.
(A proof of stronger statements can be found in \cite{BR}.)

\begin{prop}\label{propnoncommprod} Let $A$ be a rational formal power series in $k$ 
non-commuting variables. Then the shuffle product of $A$
with $\frac{1}{1-(x_0+\dots+x_k)}$
and the ordinary non-commutative product
$\frac{1}{1-(x_0+\dots+x_k)}A\frac{1}{1-(x_0+\dots+x_k)}$
are both rational.
\end{prop}

The shuffle-product counts subsequences encoded and weighted by $A$ 
in $k-$ary expansions
of natural integers
and the ordinary product in Proposition \ref{propnoncommprod} counts
subfactors (encoded by $A$) in $k-$ary expansions.

Proposition \ref{propotherdescr} corresponds to the case 
where 
$$A=x_1\frac{1}{1-x_0x_1}=x_1+x_1x_0x_1+x_1x_0x_1x_0x_1+x_1x_0x_1x_0x_1+
\dots$$ 
respectively 
$$A=x_1\frac{w}{1-x_0x_1}=x_1+wx_1x_0x_1+w^2x_1x_0x_1x_0x_1+\dots$$
in the weighted case with $x_0,x_1$ non-commuting
variables and $w$ a central variable.

A famous example counting subsequences or subfactors (reduced to $1$) 
is given by the Thue-Morse 
sequence corresponding to $A=x_1$. 

Another famous example counting subfactors is given by the 
Rudin-Shapiro sequence associated to $A=x_1^2$.


\section{Proof of Theorem \ref{thmstid} and more formulae}\label{sectproofsmain}

\subsection{Proof of Theorem \ref{thmstid}}
{\bf Proof of assertion (i)} For $n=0$ we have
$$s(1+0)=s(1)+s(0)=1+0=1,\ s(1+1)=s(0)+s(1)=0+1=1$$
and $$t(1+0)=s(1)-s(0)=1-0=1,\ t(1+1)=s(0)-s(1)=0-1=-1\ .$$
The proof is now by induction on $e$. If $n$ is even we have
$$s(2^e+n)=s(2^{e-1}+\frac{n}{2})=s(2^{e-1}-\frac{n}{2})+s(\frac{n}{2})
=s(2^e-n)+s(n)$$
and
\begin{eqnarray*}
t(2^e+n)&=&-t(2^{e-1}+\frac{n}{2})\\
&=&-(-1)^{e-1}\left(s(2^{e-1}-\frac{n}{2})-s(\frac{n}{2})\right)\\
&=&(-1)^e(s(2^e-n)-s(n))
\end{eqnarray*}

If $n$ is odd, we have
\begin{eqnarray*}
s(2^e+n)&=&s(2^{e-1}+\frac{n-1}{2})+s(2^{e-1}+\frac{n+1}{2})\\
&=&s(2^{e-1}-\frac{n-1}{2})+s(\frac{n-1}{2})+s(2^{e-1}-\frac{n+1}{2})+
s(\frac{n+1}{2})\\
&=&s(2^e-n)+s(n)
\end{eqnarray*}
and 
\begin{eqnarray*}
&&t(2^e+n)\\
&=&-t(2^{e-1}+\frac{n-1}{2})-t(2^{e-1}+\frac{n+1}{2})\\
&=&(-1)^e\left(s(2^{e-1}-\frac{n-1}{2})-s(\frac{n-1}{2})+s(2^{e-1}
-\frac{n+1}{2})-s(\frac{n+1}{2})\right)\\
&=&(-1)^e\left(s(2^e-n)-s(n)\right)
\end{eqnarray*}

{\bf Proof of assertion (ii)} These formulae are easy to establish
for $e=0$. 

For even $n$ we have
$$t(3\cdot 2^e+n)=-t(3\cdot 2^{e-1}+\frac{n}{2})=
-(-1)^{e-1}s(\frac{n}{2})=(-1)^es(n)$$
and 
$$t(6\cdot 2^e-n)=-t(6\cdot 2^{e-1}-\frac{n}{2})=
-(-1)^{e-1}s(\frac{n}{2})=(-1)^es(n)$$
and for odd $n$ we get
\begin{eqnarray*}
t(3\cdot 2^e+n)&=&-t(3\cdot 2^{e-1}+\frac{n+1}{2})-t(3\cdot 2^{e-1}+\frac{
n-1}{2})\\
&=&-(-1)^{e-1}(s(\frac{n+1}{2})+s(\frac{n-1}{2}))=(-1)^es(n)
\end{eqnarray*}
and
\begin{eqnarray*}
t(6\cdot 2^e-n)&=&-t(6\cdot 2^{e-1}-\frac{n+1}{2})-t(6\cdot 2^{e-1}-\frac{
n-1}{2})\\
&=&-(-1)^{e-1}(s(\frac{n+1}{2})+s(\frac{n-1}{2}))=(-1)^es(n)\ .
\end{eqnarray*}
This completes the proof.
\hfill$\Box$

\subsection{A few other formulae}

\begin{prop} \label{propmoreform1} (i) We have
$$s(2^{e+1}+n)=s(2^e+n)+s(n)$$
for $0\leq n\leq 2^e$, (see the remark by T. Tokita concerning 
the Stern-sequence A2487 in \cite{OEIS}). 

(ii) We have
$$t(2^{e+1}+n)+t(2^e+n)=(-1)^{e+1}s(n)$$
for $0\leq n\leq 2^e$.
\end{prop}

The formulae of Proposition \ref{propmoreform1} have the following 
conjectural generalisation, analogous to Conjecture \ref{conjgen}:

\begin{conj} (i) The series 
$$A(z)=\frac{\sum_{n=0}^\infty (s(2+n)-s(1+n))z^n}
{\sum_{n=0}^\infty s(n)z^n}=1-2z+2z^2-4z^4+4z^5+2z^6+\dots 
$$
satisfies
$$\sum_{n=0}^\infty (s(2^{e+1}+n)-s(2^e+n))z^n=A(z^{2^e})\sum_{n=0}^\infty
s(n)z^n$$
for all $e\in\mathbb N$.

Similarly, the series 
$$B(z)=-\frac{\sum_{n=0}^\infty (t(2+n)+t(1+n))z^n}
{\sum_{n=0}^\infty s(n)z^n}=1-2z-2z^2+4z^3+6z^6-6z^7+\dots 
$$
satisfies 
$$(-1)^{e+1}\sum_{n=0}^\infty (t(2^{e+1}+n)+ t(2^e+n))z^n=B(z^{2^e})\sum_{n=0}^\infty
s(n)z^n$$
for all $e\in\mathbb N$.
\end{conj}

{\bf Proof of Proposition \ref{propmoreform1}} 
The formulae hold for $n=0$ and $a\in\{0,1\}$.
The induction step is an easy computation for odd $a$ and obvious for even
$a$.\hfill$\Box$

\begin{prop} (i)
We have 
$$s(n)=-s(n-2^e)+s(n-2\cdot 2^e)+2s(n-3\cdot 2^e),\ 2^{e+2}\leq n\leq
2^{e+3}-2^e$$
for $e\geq 0$.

(ii) We have
$$t(n)=t(n-2^e)-t(n-2^{e+1})$$
for $2^{e+2}\leq n\leq 2^{e+3}$.
\end{prop}

{\bf Proof} The case $e=0$ implies $n\in\{4,5,6,7\}$
in assertion (i) and we have
$$\begin{array}{l}
s(4)=1=-2+1+2\cdot 1=-s(3)+s(2)+2s(1)\\
s(5)=3=-1+2+2\cdot 1=-s(4)+s(3)+2s(2)\\
s(6)=2=-3+1+2\cdot 2=-s(5)+s(4)+2s(3)\\
s(7)=3=-2+3+2\cdot 1=-s(6)+s(5)+2s(4)\end{array}$$

The induction step for $e>0$ is easy if $n$ is even and 
involves the usual identity $s(n)=s((n-1)/2)+s((n+1)/2)$ if
$n$ is odd.

The proof of assertion (ii) is similar.\hfill$\Box$

The following result gives a few partial sums associated to
the Stern sequence and its twist:

\begin{prop} We have
\begin{eqnarray*}
\sum_{n=1}^{2^e}s(n)&=&\frac{3^e+1}{2},\ e\geq 0\\
\sum_{n=1}^{2^e}(-1)^ns(n)&=&\frac{1-3^{e-1}}{2},\ e\geq 1\\
\sum_{n=1}^{2^e}t(n)&=&\frac{(-1)^e+1}{2},\ e\geq 0\\
\sum_{n=1}^{2^e}(-1)^nt(n)&=&\frac{-3+(-1)^e}{2},\ e\geq 0\ .
\end{eqnarray*}
\end{prop}

{\bf Proof} The first equality holds for $e=0$ and 
\begin{eqnarray*}
\sum_{n=0}^{2^{e+1}}s(n)&=&\sum_{n=0}^{2^e}s(2n)+\sum_{n=1}^{2^e}
s(2n-1)\\
&=&\sum_{n=0}^{2^e}s(n)+\sum_{n=1}^{2^e} (s(n-1)+s(n))\\
&=&-s(2^e)+3\sum_{n=0}^{2^e}s(n)=-1+3\frac{3^e+1}{2}=
\frac{3^{e+1}+1}{2}
\end{eqnarray*}
by induction.

For the next identity one finds similarly
$$\sum_{n=1}^{2^{e+1}}(-1)^ns(n)=1-\sum_{n=1}^{2^e}s(n)=1-\frac{3^e+1}{2}=
\frac{1-3^e}{2}\ .$$

The computations for the partial sums involving
$t(n)$ and $(-1)^n t(n)$ are analogous.
\hfill$\Box$

We end this section with a list of a few more identities.

\begin{prop} We have
\begin{eqnarray}
s(3\cdot 2^e+n)=s(3\cdot 2^e-n)
\end{eqnarray}
for all $e,n$ such that $0\leq e\leq 2^n$,
\begin{eqnarray}
s(3\cdot 2^e+n)=s(3\cdot 2^{e-1}+n)+2s(n)
\end{eqnarray}
for all $e,n$ such that $0\leq n\leq 2^{e-1}$,
\begin{eqnarray}
t(2^e+n)=t(2^e+n-2^{e-2})-t(2^e+n-2^{e-1})
\end{eqnarray}
for all $e,n$ such that $e\geq 2$ and $1\leq n\leq 2^e$,
\begin{eqnarray}
s(2^e+n)=(-1)^et(2^e+n)+2s(n)\end{eqnarray}
for all $e,n$ such that $0\leq n\leq 2^{e+1}$,
\begin{eqnarray}
s(2^e+n)=(-1)^et(2^e-n)-3s(n)
\end{eqnarray}
for all $e,n$ such that $0\leq n\leq 2^{e-1}$,
\begin{eqnarray}
s(2^e-n)=(-1)^et(2^e-n)+2s(n)
\end{eqnarray}
for all $e,n$ such that $0\leq n\leq 2^{e-1}$,
\begin{eqnarray}
s(2^e-n)=(-1)^et(2^e+n)+s(n)
\end{eqnarray}
for all $e,n$ such that $0\leq n\leq 2^{e}$.
\end{prop}

Proofs are easy and left to the reader.


\section{Proofs related to factorisations}\label{sectproofkreg}

{\bf Proof of Theorem \ref{thmtpoly}} We set 
$\psi_e=z(1+z^{2^e})(1+z+z^2)\prod_{n=0}^{e-2}(1-z^{2^n}
+z^{2^{n+1}})^{e-1-n}$. Iterating the trivial identity 
$$(1+z^n+z^{2n})(1-z^n+z^{2n})=
(1+z^{2n}+z^{4n})$$ we get the equivalent expression
$$\psi_e=z(1+z^{2^e})\prod_{n=0}^{e-1}(1+z^{2^n}+z^{2^{n+1}})\ .$$
The proof of the identity $\psi_e=(-1)^e\sum_{n=0}^{3\cdot 2^e}
t(3\cdot 2^e+n)z^n$ is by induction on $e$. It holds for $e=0$. 
The induction step follows from the 
recursive definition of the sequence $t(0),t(1),\dots$ and from
the equality $\psi_{e+1}(z)=\left(\frac{1}{z}+1+z\right)\psi_e(z^2)$.
\hfill$\Box$

{\bf Proof of Theorem \ref{thmlogder}:}
Using the Carlitz factorisation 
$$\tilde S(z)=\prod_{n=0}^\infty \left(1+z^{2^n}+z^{2^{n+1}}\right)$$
of $\tilde S(z)=\sum_{n=0}^\infty s(n+1)z^n$ we have
$$H(z)=\frac{d}{dz}\mathop{log}(\tilde S(z))=
\sum_{n=0}^\infty\frac{2^nz^{2^n-1}+2^{n+1}z^{
2^{n+1}-1}}{1+z^{2^n}+z^{2^{n+1}}}\ .$$
The summand of index $n=0$ yields $\frac{1+2z}{1+z+z^2}$ and the sum 
$\sum_{n=1}^\infty\cdots$ can be rewritten as $2zH(z^2)$.

The proof of $2-$regularity of $H(z)$ is an easy consequence of the functional
equation for $H$, see Theorem \ref{thmaffpres} below.

Uniqueness of $H$ defined by the functional equation
$H(z)=\frac{1+2z}{1+z+z^2}+2zH(z^2)$ follows from the fact that
the map 
$$A(z)\longmapsto \frac{1+2z}{1+z+z^2}+2zA(z^2)$$
has a unique attracting fixpoint for formal power series (with
respect to the obvious topology given by coefficent-wise 
convergency).
\hfill$\Box$

{\bf Proof of Theorem \ref{thmaffpres}} We assume first that 
no linear form $L_1,\dots,L_d$ involves coefficients of degree
$\geq k$.
For $i=0,\dots,k-1$, we denote as before
by $\rho(i)$ the linear map
$$\rho(i)\left(\sum_{n=0}^\infty a(n)z^n\right)=\sum_{n=0}^\infty
a(i+nk)z^n\ .$$
Given solutions $U_1,\dots,U_d$, we consider a finitely generated
vector space or module $\mathcal V$ containing $U_1,\dots,U_d$ and the 
$k-$kernel of $A_1,\dots,A_d$. We have then
$$\rho(i)U_j=\rho(i)A_j+([x^i]L_j)(U_1,\dots,U_d)$$
where $[x^i]L_j\in R[x_1,\dots,x_d]$ is the linear form obtained from 
$L_j$ by considering the coefficients of $z^j$.
The power series $\rho(i)U_j\in\mathcal V$ is thus a linear 
combination of $U_1,\dots,U_1$ and of the $k-$kernel of $A_j$.
The set $\mathcal V$ contains thus the $k-$kernel of $U_j$ and
$U_1,\dots,U_d$ are all $k-$regular. 

If there are linear forms among $L_1,\dots,L_d$ which are of degree
$\geq k$, we introduce the $k-$regular series $A_{d+1}=zA_1,\dots,
A_{2d}=zA_d$, the series $U_{d+1}=zU_1,\dots,U_{2d}=zU_d$,
the linear forms $L_{d+1}=zL_1,\dots,L_{2d}=zL_d\in zR[z][x_1,\dots,x_d]$
and set $\alpha_{d+1}=\dots=\alpha_{2d}=0$. We have then
the equations 
$$U_i(z)=A_i(z)+zL_i(U_1(z^k),\dots,U_n(z^k))$$
and the identities $\alpha_i=U_i(0)$
for $i=1,\dots,2d$. Modifying a linear form $L_j$ of degree $\geq k$ 
by substituting all occurences of $z^{k+i}U_j(z^k)$ with
$z^iU_{d+j}(z^k)$ we construct an equivalent system 
with strictly smaller maximal degree for the linear forms $L_1,\dots,L_{2d}$.
Iteration of this construction leads eventually to a system
containing only linear forms of degree strictly smaller than $k$. 

Existence and unicity of the solution follow from unicity of 
the attracting fixpoint of the dynamical system defined by the map
$$U_i(z)\longmapsto A_i(z)+L_i(U_1(z^k),\dots,U_d(z^k)),i=1,\dots,d$$
starting from the point $(\alpha_1,\dots,\alpha_d)$.
\hfill$\Box$

{\bf Proof of Theorem \ref{thmPreg}} 
The first part follows from Theorem \ref{thmaffpres} applied to the 
identity $A(z)=P(z)A(z^k)$. 
We present here however a second, independent proof.

Working over the field of 
complex numbers and using the fact that products of two
$k-$regular series are $k-$regular (cf. Theorem 16.4.1 in \cite{AS}), 
it is enough to prove the result for polynomials of degree
$1$. We can thus assume that $P(z)=1+\lambda z$. 
The coefficient of $z^n$ in 
$A(z)=\prod_{m=0}^\infty (1+\lambda z^{k^m})$ is then given by zero
if the $k-$ary expansion of $n$ involves digits greater than $1$ and 
it is given by $\lambda^\alpha$ otherwise, 
where $\alpha$ equals the number of ones in the $k-$ary expansion of $n$.
This implies $k-$regularity of $A(z)$.

We have 
$$B(z)=\frac{d}{dz}\mathop{log}(A(z))=\sum_{n=0}\frac{P'(z^{k^n})k^nz^{
k^n-1}}{P(z^{k^n})}\ .$$
The summand of index $n=0$ yields $\frac{P'(z)}{P(z)}$ and the remaining
summation $\sum_{n=1}^\infty\cdots$ can be rewritten as $kz^{k-1}B(z^k)$.
This shows that $B(z)$ satisfies the functional
equation
$$B(z)=\frac{P'(z)}{P(z)}+kz^{k-1}B(z^k)\ .$$
Using Theorem 16.4.3 of \cite{AS} we see that the rational 
fraction $\frac{P'(z)}{P(z)}$ 
is $k-$regular if and only if all zeroes of $P(z)$ are roots of unity
(ie. if $P(z)$ divides $(z^N-1)^N$
for some integer $N$).
Theorem \ref{thmaffpres} implies then $k-$regularity of $B(z)$.
\hfill$\Box$


\section{Proof of Theorem \ref{thmmatr} and other results involving matrices}
\label{sectproofmatr}

{\bf Proof of Theorem \ref{thmmatr}}
The trivial identities
\begin{eqnarray*}
\det(M(2n))&=&\det\left(\begin{array}{cc}
s(2n)&s(2n+1)\\t(2n)&t(2n+1)\end{array}\right)\\
&=&\det\left(\begin{array}{cc}
s(n)&s(n)+s(n+1)\\-t(n)&-t(n)-t(n+1)\end{array}\right)\\
&=&-\det\left(\begin{array}{cc}
s(n)&s(n+1)\\t(n)&t(n+1)\end{array}\right)\\
&=&-\det(M(n))
\end{eqnarray*}
and 
\begin{eqnarray*}
\det(M(2n-1))&=&\det\left(\begin{array}{cc}
s(2n-1)&s(2n)\\t(2n-1)&t(2n)\end{array}\right)\\
&=&\det\left(\begin{array}{cc}
s(n-1)+s(n)&s(n)\\-t(n-1)-t(n)&-t(n)\end{array}\right)\\
&=&-\det\left(\begin{array}{cc}
s(n-1)&s(n)\\t(n-1)&t(n)\end{array}\right)\\
&=&-\det(M(n-1))
\end{eqnarray*}
imply the result.\hfill$\Box$

{\bf Proof of Proposition \ref{propstmod2}} 
The reduction modulo $2$ of the Stern sequence $s(0),s(1),\dots$
is the $3-$periodic sequence $0,1,1,0,1,1,\dots$. 
Indeed, this holds for $s_0=0,s_1=s_2=1$
and the recursive formulae
\begin{eqnarray*}
s(6n)&=&s(3n)\\
s(6n+1)&=&s(3n)+s(3n+1)\\
s(6n+2)&=&s(3n+1)\\
s(6n+3)&=&s(3n+1)+s(3n+2)\\
s(6n+4)&=&s(3n+2)\\
s(6n+5)&=&s(3n+2)+s(3n+3)
\end{eqnarray*}
imply the $3-$periodicity of $s(n)\pmod 2$ by induction.
The reduction modulo $2$ of twisted Stern sequence $t(0),t(1),\dots$
coincides with the reduction modulo $2$ of the Stern sequence.
\hfill$\Box$

\subsection{Other results involving matrices}

The proofs of the following results are easy and omitted.

\begin{prop}
(i) The matrices 
$$\left(\begin{array}{cc}s(n)&s(n+1)\\s(2^e+n)&s(2^e+n+1)\end{array}
\right)$$
have determinant $-1$ for $n$ such that $0\leq n<2^e$ and determinant $1$ for
$n$ such that $2^e\leq n<2^{e+1}$.

\ \ (ii) The matrices 
$$
\left(\begin{array}{cc}s(n)&s(n+1)\\t(2^e+n)&t(2^e+n+1)\end{array}
\right)$$
have determinant $(-1)^{e+1}$ for $n$ such that $0\leq n<2^e$ and determinant
$(-1)^e$ for $n$ such that $2^e\leq n<2^{e+2}$.

\ \ (iii) The matrices 
$$\left(\begin{array}{cc}t(n)&t(n+1)\\s(2^e+n)&s(2^e+n+1)\end{array}
\right)$$
have determinant $(-1)^{e+1}$ for $n$ such that $2^{e+1}< n<5\cdot 2^e$.

\ \ (iv) The matrices 
$$\left(\begin{array}{cc}t(n)&t(n+1)\\t(2^e+n)&t(2^e+n+1)\end{array}
\right)$$
have determinant $1$ for $n$ such that $2^{e-2}\leq n<2^e$ or 
$7\cdot 2^e\leq n <2^{e+3}$
and determinant
$-1$ for $n$ such that $2^e\leq n<7\cdot 2^e$.
\end{prop}


\section{Proof for Theorem \ref{thmdiv}}\label{sectdiv}

For odd $n$ we have
$$s(n)=s((n-1)/2)+s((n+1)/2)=s(n-1)+s(n+1)\ .$$
We have thus $\frac{s(n-1)+s(n+1)}{s(n)}=1=1+2v_2(n)$
since $v_2(n)=0$ if $n$ is odd.

For $n$ even we have by induction
\begin{eqnarray*}
&&\frac{s(n-1)+s(n+1)}{s(n)}\\
&=&\frac{s((n-2)/2)+s(n/2)+s(n/2)+s((n+2)/2)}{s(n/2)}\\
&=&\frac{s(n/2-1)+s(n/2+1)}{s(n/2)}+2\\
&=&1+2v_2(n/2)+2=1+2v_2(n)\ .
\end{eqnarray*}
This ends the proof of assertion (i).

For the twisted Stern sequence we use the analogous identities
\begin{eqnarray*}
t(n-1)+t(n+1)&=&t(n),\ n\hbox{ odd,}\\
t(n-1)+t(n+1)&=&-(t(n/2-1)+t(n/2+1)+2t(n)),\ n\hbox{ even}, n\geq 4.
\end{eqnarray*}
This implies assertion (ii) by checking the initial cases
and the case of $n\in 3\cdot 2^{\mathbb N}$.
\hfill$\Box$

\noindent Roland BACHER

\noindent INSTITUT FOURIER

\noindent Laboratoire de Math\'ematiques

\noindent UMR 5582 (UJF-CNRS)

\noindent BP 74

\noindent 38402 St Martin d'H\`eres Cedex (France)
\medskip

\noindent e-mail: Roland.Bacher@ujf-grenoble.fr


\begin{thebibliography}{99}

\bibitem{AZ} M. Aigner and G. M. Ziegler, Proofs from THE BOOK, 
3rd ed., Springer-Verlag (2004).

\bibitem{AS} J.-P. Allouche, J. Shallit, Automatic Sequences.
Theory, Applications, Generalizations, Cambridge
University Press (2003).

\bibitem{BR} J. Berstel and C. Reutenauer, {\it Noncommutative Rational 
Series with Applications}, available at the authors websites. 

\bibitem{dB} G. de Bruijn, {\it On Mahler's partition problem}, Indag. Math. 
vol. 10 (1948), 210--220.

\bibitem{CW} N. Calkin, H. S. Wilf, {\it Recounting the rationals}, Amer. Math. Monthly, {\bf 107} (2000), 360--363.

\bibitem{C} L. Carlitz, {\it A problem in partitions related to the 
Stirling numbers}, Bull. Amer. Math. Soc., 70(2) (1964), 275--278.
2319565 (2008c:11033)

\bibitem{KMP} S. Klavzar, U. Milutinovic, C. Petr, {\it 
Stern polynomials},
Adv. in Appl. Math. 39 (2007), no. 1, 86--95.

\bibitem{L} D. H. Lehmer, {\it On Stern's Diatomic Series}, 
Amer. Math. Monthly 36(1) 
1929, 59--67.

\bibitem{OEIS} N. J. A. Sloane, (2008), The On-Line Encyclopedia of 
Integer Sequences,
published electronically at www.research.att.com/~njas/sequences/.

\bibitem{S} M. A. Stern, {\it \"Uber eine zahlentheoretische Funktion}, 
J. Reine Angew. Math., 55 (1858), 193--220. 

\end{thebibliography}
\end{document}